\documentclass[notitlepage,leqno]{article}

\usepackage{latexsym}
\usepackage{amsmath}
\usepackage{amsfonts}
\usepackage{amssymb}
\renewcommand{\a }{\alpha }

\renewcommand{\d}{\delta }
\newcommand{\D }{\Delta }

\newcommand{\e }{\varepsilon }
\newcommand{\g }{\gamma}

\renewcommand{\l }{\lambda }

\newcommand{\n }{\nabla }
\newcommand{\var }{\varphi }

\newcommand{\s }{\sigma }
\newcommand{\Sig }{\Sigma}

\renewcommand{\th }{\theta }
\renewcommand{\o }{\omega }

\newcommand{\ov}{\overline}

\newcommand{\be}{\begin{equation}}
\newcommand{\ee}{\end{equation}}
\newenvironment{pf}{\noindent{\sc Proof}.\enspace}{\rule{2mm}{2mm}\medskip}
\newenvironment{pfn}{\noindent{\sc Proof}}{\rule{2mm}{2mm}\medskip}

\newcommand{\R}{\mathbb{R}}
 \newcommand{\salt}{\noalign{\vskip .2truecm}}
 \newcommand{\parent}[3]{\left #1 {#3} \right #2} % racchiude il testo tra il
\newcommand{\barre}[1]{\parent \Vert \Vert {#1}} % racchiude il testo tra barre adeguate

\newcommand{\Z}{\mathbb{Z}}

\newcommand{\N}{\mathbb{N}}

\newtheorem{pro}{Proposition}[section]
\newtheorem{thm}{Theorem}[section]
\newtheorem{lem}{Lemma}[section]

\newtheorem{df}{Definition}[section]

\author{{\sc  Mohameden Ould Ahmedou}}
\title{ \Large \textbf{A Riemann mapping type Theorem in higher dimensions  \\
Part I : the conformally flat case with umbilic boundary}}

\begin{document}

\date{ }

\maketitle

{\footnotesize
\begin{abstract}

\noindent
In this paper we prove that every Riemannian metric on a locally conformally flat manifold with umbilic boundary can be conformally deformed to a scalr flat metric having constant mean curvature. This result can be seen as a generalization to higher dimensions of the well known Riemann mapping Theorem in the plane.

\bigskip\bigskip

\end{abstract}
}

\section{Introduction}

In \cite{E1}, Jos\'e F. Escobar raised the following question:
Given a compact Riemannian manifold with boundary, when it is conformally equivalent to one that has zero scalar curvature and whose boundary has a constant mean curvature ?
This problem can be seen as a ``generalization'' to higher dimensions of the well known Riemannian mapping Theorem. The later states that an open, simply connected proper subset of the plane is conformally diffeomorphic to the disk. In higher dimensions few regions are conformally diffeomorphic to the ball. However one can still ask whether a domain is conformal to a manifold that resembles  the ball into ways : namely, it has zero scalar curvature and its boundary has constant mean curvature. In the above the term ``generalization'' has to be understood in that sens.
The above problem is equivalent to finding a smooth positive solution to the following nonlinear boundary value problem on a Riemannian manifold with boundary $ (M^n,g)$, $ n \geq 3$:  
\begin{equation}\tag{$P$}
\begin{cases}
- \Delta_g u + \frac{(n-2)}{4(n-1)}R_g u = 0, \quad u > 0  & 
\mbox{ in } \mathring{M}; \\
\partial_\nu u + \frac{n - 2}{2} h_g u = Q(M, \partial M) u^{\frac{n}{n-2}}, & 
\mbox{ on } \partial M.  
\end{cases}
\end{equation}
 where $R$ is the scalar curvature of $M$, $h$ is the mean curvature of $ \partial M$, $ \nu$ is the outer normal vector with respect to $g$ and  $ Q(M,\partial M)$ is a constant whose sign is uniquely determined by the conformal structure.
Indeed  if $ \ov{g} = u^{\frac{4}{n - 2}} g$, then the  metric $ \ov{g}$ has zero scalar curvature and the boundary has constant mean curvature with respect to $ \ov{g}$.

Solutions of equation (P) correspond , up to  a multiple constant, to critical points of the following functional $ J$ defined on $ H^1(M) \setminus \{0 \}$
\begin{equation}\label{eq:J}
 J(u ) = \, \frac{ \left(\int_M \left( |\n_g u|^2 + \frac{n-2}{4(n-1)}
\, R_g \, u^2 \right)dV_g \,  + \frac{n-2}{2} \int_{\partial M} h_g \, u^2 d\s_g \right)^{\frac{n - 1}{n - 2}}
}{  \int_{\partial M}   \,
|u|^{2\frac{n-1}{n-2}} d\s_g}. 
\end{equation}
where $ dV_g$ and $ d \s_g$ denote the Riemannian measure on $M$ and $ \partial M$ induced by the metric $g$.

The regularity of the $ H^1$ solutions of (P) was established by P. Cherrier \cite{ch}; and related problems regarding conformal deformations of metrics on manifold with boundary were studied  in \cite{alm} , \cite{[ACPY3]} , \cite{DMA}, \cite{VA} , \cite{ha}, \cite{HL1} , \cite{HL2}, \cite{ll}  , \cite{ma} and the references therein.

The exponant $\frac{2(n - 1)}{n - 2}$ is critical for the Sobolev trace embedding $ H^1(M) \to L^{q}(\partial M)$. This embedding being not compact , the functional $J$ does not satisfy the Palais Smale condition. For this reason standard variational methods cannot be applied to find critical points of $J$.

Following the original arguments introduced by T. Aubin \cite{auw}, \cite{aul} and R. Schoen \cite{Sc} to prove Yamabe conjecture on closed manifolds, Escobar proved the existence of a smooth positive solution $u$ of (P) on $ (M^n, g) , n \geq 3$ for many cases.
To state his results we need some preliminaries:

Let $ H$ denote the second fondamental form of $ \partial M$ in $(M,g) $ with respect to the inner normal. Let us denote the traceless part of the second fundamental form by $U $ that is  $U(X,Y) = \, H(X,Y) \, - \, h_g \, g(X,Y) $

\begin{df}
A point $q \in \partial M$ is called an  { \bf  umbilic point }  if $ U = 0$ at $ q$.
$ \partial M$ is called umbilic if every point of $\partial M $ is umbilic.
\end{df}
 Regarding the above problem Escobar  proved the following Theorem \cite{E1,E2}: 
\begin{thm}
Let $(M^n,g)$ be a compact Riemannian manifold with boundary, $ n \geq 3$. 
Assume that $ M^n$ satisfies one of the following conditions:
\begin{description}
\item(i) $ n \geq 6$ and $M$ has a nonumbilic point on $ \partial M$
\item(ii)  $ n \geq 6$ and $M$ is conformally locally flat with umbilic boundary
\item(iii) $ n = 4,5$ and $\partial M$ is umbilic
\item(vi) $ n = 3$
\end{description}
then there exists a smooth metric $ u^{\frac{4}{n - 2}} \, g$, $ u > 0 \,  \mbox{on } M$ of zero scalar curvature and constant mean curvature on $ \partial M$.
\end{thm}

In his proof Escobar uses strongly an extension of the positive mass Theorem of R. Schoen and S.T. Yau  \cite{sy}, \cite{sz} to some type of manifolds with boundary. Such an extension was proved by Escobar in \cite{E4}.
Besides the proof of T.Aubin and R.Schoen of the Yamabe conjecture, another proof by A. Bahri \cite{ba} and A.Bahri and H. Brezis \cite{bb} of the same conjecture is available by techniques related to the Theory of critical point at Infinity of A. Bahri \cite{bal}.

We plan to give a complete positive answer to the above problem based on the topological argument of Bahri-Coron \cite{bc}, as Bahri and Brezis did for the Yamabe conjecture. In this first part we study the case where the manifold is locally conformally flat with umbilic boundary. Namely we prove the following Theorem 

\begin{thm}\label{t:mt}
Suppose that $(M^n,g)$, $ n \geq 3$ is a compact locally conformally flat manifold with umbilic boundary, then equation (P) has a solution.
\end{thm}
 Let us observe that while the solution obtained by Escobar is a minimum of $J$,  our solution is in general, a critical point of $J$ of higher Morse index, more precisely we have the following characterization of  solutions obtained  by Bahri-Coron existence scheme( see  \cite{yc}) :
\begin{thm}\label{t:ct}
The solution  $u$ obtained in Theorem \ref{t:mt} satisfies, for some nonnegative integer $p_0$:
\begin{description}
\item(i) $ p_0^{\frac{1}{n - 2}} S\leq J(u) \leq (p_0 + 1)^{\frac{1}{n - 2}}S$
\item(ii) $ind(J,u) \leq (p_0 + 1)( n - 1) + p_0, \, ind(J,u) + dim ker d^2 \, J(u) \geq p_0(n - 1  ) + p_0 $
\item(iii) $ u$ induces some difference of topology at the level $ p_0^{\frac{1}{n - 2}} S$.
\end{description}
where $ind(J,u)$ is the Morse index of $J$ at $u$, and $ S= 2^{1 - n} \o_{n - 1}$, where $ \o_{n - 1}$ is the volume of the $n - 1 $ dimensional unit sphere.
Moreover, if (P ) has only nondegenerate solutions, $ ind(J,u) = (n - 1)p_0 + p_0 $.
\end{thm}

The remainder of the paper is organized as follows: in section 2 we construct some ``almost solutions'' which are solutions of the ``problem at Infinity''. In section 3 we collect some standard  results regarding the description of the lack of compactness and some local deformation Lemma. In section 4 we perform an expansion of $J$ at 'Infinity' and we give the proof of Theorem \ref{t:mt} in section 5. Lastly we devote the appendix to establish some technical Lemma and to recall some well known results.

\begin{center}

{\bf Acknowledgements}

\end{center}

The author is indebted to Pr. Abbas Bahri for teaching him his Theory of critical point at Infinity and he is gratefull to Pr. Antonio Ambrosetti for his interest in his work and his constant support.

\section{Construction of ``almost solutions''}

In this paper we assume that $(M^n,g)$ is a compact Riemannian manifold with boundary and dimension $ n \geq 3 $ . Let $ R_{pq}$ and $R = g^{pq}R_{pq} $ be the Ricci curvature and the scalar curvature, respectively; let $ h_{ij}$ and $ h = \frac{1}{n - 1}g^{ij}\, h_{ij}$ be the second fundamental form of the boundary of M, $\partial M$ and the mean curvature , respectively. Let $ \tilde{g} = u^{\frac{4}{n - 2}} g$ be a metric conformally related to $ g$.
We denote by a tilde all quantities computed with respect to the metric $ \tilde{g}$. The transformation law for the scalar curvature is 
\begin{equation}\label{e:r}
\tilde{R} = \frac{4(n - 1)}{n - 2} \, \frac{L u}{u^{\frac{n + 2}{n - 2}}}
\end{equation}
where $ L$ is the conformal Laplacian $ L = \D \, - \, \frac{n - 2}{4(n - 1)} R \quad \mbox{on } M $; while the transformation law for the mean curvature is
\begin{equation}\label{e:h}
\tilde{h} = \frac{2}{n - 2} \, \frac{B u}{u^{\frac{n}{n - 2}}}
\end{equation}
where $B $ is the boundary operator $  B = \frac{\partial}{\partial \nu} \, + \, \frac{n - 2}{2} h \quad \mbox{on } \partial M$.

Consider now the following eigenvalue problem on $(M,g) $ :
\begin{equation}\tag{$E$}
\begin{cases}
L_g u = \l \,  u  \quad \mbox{on } \, \mathring{M} \\
B_g u = 0 \quad \mbox{on } \, \partial M
\end{cases}
\end{equation}
Let $ \l_1 $ the first eigenvalue of (E).

\begin{df}
We say that a manifold $ M$ is of positive(negative, zero) type if $ \l_1 > 0 (< 0, = 0 )$.
\end{df}

As it is well known the  existence problem is easy when the manifold is of negative or zero type, so we treat only in this paper the case of manifold of positive type.

\noindent Now we construct some almost solutions of (P), which will play a central role in the description of the lack of compactness.

Let $f_1$ denote a positive eigenfunction corresponding to the first eigenvalue of (E) , and consider $ g_1 = (f_1)^{\frac{4}{n - 2}} g$ 
 then , according to ($ \ref{e:r}$) and ($\ref{e:h}$) we have: $ R_{g_{1}} > 0 $ and $ h_{g_1} = 0$ on $ \partial M$. We can work with $g_1$ instead of $g$, but for simplicity we still denote it by $g$.
Let $a \in \partial M$; since $M$ is a compact  locally conformally flat manifold one can find a neighborhood of $a$,  $\mathcal{U}(a) \supset  B^M_{\rho}(a)$  , $\rho > 0$ uniform and a conformal diffeomorphism $\varphi$ which maps $ B^M_{\rho}(a)$ into $\R^n$ with $\varphi(0) = a $ . Therefore,  denoting $g_0$ the flat metric on $\R^n$, there exist a positive function $u_a$ such that $ \varphi^{*}(g_0) = u_a^{\frac{4}{n - 2}} \, g$.
Since the boundary is umbilic, $\varphi (\partial M \cap B^M_{\rho}(a))$ has to be a piece of sphere or a piece of a hyperplane (See \cite{sp}) and since spheres and hyperplanes are locally conformal to each other, we can assume without loss of generality that $ \partial B^+_2(0) \cap \partial \R^n_+ \subset \varphi(\partial M \cap
B^M_{\rho}(a) )$ and $ \varphi(\mathring{ M} \cap
B^M_{\rho}(a) ) \subset \R^n_+$ .
 Since $\partial B^+_2(0) \cap \partial \R^n_+ $ has zero mean curvature in $ \ov{B^+_2}$, we deduce from  $(\ref{e:h})$ that $  \frac{\partial u_a}{\partial \nu}  = 0$ on $\partial M \cap B^M_{\rho}(a) $. We extend $ u_a$ to be a smooth positive function on $M$ such that $  \frac{\partial u_a}{\partial \nu}  = 0$ on $ \partial M$ and $u_a = 0 $ on $ M \setminus B^M_{2\rho}(a) $.
Consider now the conformal metric $ \ov{g_0} = u^{\frac{4}{n - 2}} \, g$, then $ \ov{g_0}$ has the property that $ h_{\ov{g_0}} = 0 $ and it is Euclidean in $ B^M_{\rho}(a)$. Moreover this metric can be chosen to depend smoothly on $ a$ (see \cite{ba}).

For $a \in \partial M$, define the function:
$$
\d_{a, \l}\, (y) = \ov{c} \, \frac{\l^{\frac{n - 2}{2}}}{((1 + \l x^n)^2 + \l^2|x'|^2)^{\frac{n - 2}{2}}}
$$
where $ (x',x^n) = \varphi(y)$, and $\ov{c}$ is chosen such that $\d_{a,\l}$ satisfies the following equation
$$
\begin{cases}
- \Delta_{\ov{g}_0} u = 0, & \mbox{in } B^M_{\rho} \cap \mathring{M}; \\
\partial_\nu u  =   \, \d_{a,\l}^{\frac{n}{n - 2}} ,& \mbox{on } B^M_{\rho} \cap \partial M
\end{cases}
$$
Set $ \hat{\d}_{a, \l} = \o_a \, u_a \, \d_{a, \l}$ where $ \o_a$ is a cutoff function $ \o_a = 1 \quad \mbox{on } B^M_{\rho}(a)$ and $ \o_a = 0 \quad \mbox{on } M \setminus B^M_{2\rho}  $.

We  define now  a familly of almost solutions  $\varphi_{a,\l}$ to be the unique solution of
$$
\begin{cases}
- L_g u = 0, & \mbox{in }  \mathring{M}; \\
B_g u  =   \, \hat{\d}_{a,\l}^{\frac{n}{n - 2}} ,& \mbox{on } \partial M
\end{cases}
$$

Let us recall that the operators $ L_g$ and $ B_g$ are conformally invariant under the conformal change of metrics, namely we have:

\begin{lem}\label{l:conf}\cite{E1}

Let $\psi \in C^2(B_{\rho}(a))$, we have
$$
L_g (u_a \psi) = u_a^{\frac{n + 2}{ n - 2}}L_{\ov{g}_0}(\psi)
$$
and
$$
B_g (u_a \psi) = u_a^{\frac{n }{ n - 2}}B_{\ov{g}_0}(\psi
$$
\end{lem}

In the remainder of this section we establish some properties of our almost solutions $ \var_{a,\l}$.

\begin{lem}\label{l:error}

There are two positive constants $C$ and $B$ , such that for all $ a \in \partial M$ and $ \l \geq B$, we have
$$
{\left|\varphi_{a,\l} - \hat{\d}_{a, \l} \right|}_{\infty} \leq \frac{C}{\l^{\frac{n - 2}{2}}}
$$
\end{lem}

\begin{pf}

Let $ H_{a,\l} = \l^{\frac{n - 2}{2}} (\varphi_{a,\l} - \hat{\d}_{a,\l})$, we have
$$
L_g \, H_{a,\l} = \l^{\frac{n - 2}{2}} \quad L_g \, (\o_a u_a \d_{a,\l}) \\
= \l^{\frac{n - 2}{2}} \, u_a^{\frac{n + 2}{ n - 2}}\quad L_g \, (\o_a  \d_{a,\l})
$$
Since on $B_{\rho}$ , $\o_a = 0$, we deduce that on  $B_{\rho}$  we have $L_g \, H_{a,\l} = 0$, whereas on $M \setminus B_{\rho}$ there holds $L_g \, H_{a,\l} \leq C $.

From another part
$$
B_g \, H_{a,\l}  = \l^{\frac{n - 2}{2}} \quad [ B_g \, \varphi_{a,\l} - B_g \, (\o_a u_a \d_{a,\l})] \\
= \l^{\frac{n - 2}{2}} \,[ \hat{\d}_{a,\l} \quad - u_a^{\frac{n }{ n - 2}}\quad B_g \, (\o_a  \d_{a,\l})]
$$
on $ B_{\rho}(a) \cap \partial M, \,  \o_a = 1$, therefore $B_g \, H_{a,\l} = 0 $ , while on $M \setminus B_{\rho}$ there holds $B_g \, H_{a,\l} \leq C $.
Thus our Lemma follows from Lemma \ref{l:cp} quoted in the appendix.
\end{pf}

\begin{lem}\label{l:lb}
There are two positive constants $C$ and $B$ , such that for all $ a \in \partial M$ and $ \l \geq B$, we have
$$
\varphi_{a,\l} \,  \geq  \, \frac{C}{\l^{\frac{n - 2}{2}}}
$$
\end{lem}

\begin{pf}

Using Lemma \ref{l:error}, we know that if $ \rho_1 < \rho$ is chosen small enough, independent of $ \l$, the following inequality holds on $ B(a,\rho_1)$
$$
\varphi_{a,\l} \geq \hat{\d}_{a,\l} - \frac{C}{\l^{\frac{n - 2}{2}}} \geq \frac{C}{\l^{\frac{n - 2}{2}}}
$$ 

Let $ \Sigma_1 = \partial B(a,\rho) \cap \mathring{M}$ and $ \Sigma_2 = \partial M \setminus \Sigma_1 $.

Then we have 
\begin{equation}\label{tata}
\begin{cases}
L_g(\varphi_{a, \l} - \frac{C}{\l^{\frac{n - 2}{2}}}) \quad  \leq \, 0 & \mbox{in } \mathring{M} \\
\varphi_{a, \l} - \frac{C}{\l^{\frac{n - 2}{2}}} \quad \geq 0, & \mbox{on } \Sigma_1\\
\frac{\partial}{\partial \nu}(\varphi_{a, \l} - \frac{C}{\l^{\frac{n - 2}{2}}}) \quad \geq 0  & \mbox{on } \Sigma_2
\end{cases}
\end{equation}

Then by the hopf maximum principle, we deduce from (\ref{tata}) that
$$
\varphi_{a,\l} \geq \frac{C}{\l^{\frac{n - 2}{2}}} \quad \mbox{for } x \in M
$$

\end{pf}

\begin{lem}\label{l:fe}

Let $\theta > 0$ be given. There are positive constants $C$ and $B$, such that the following estimates hold, provided  $ \l \geq B$
\begin{description}
\item(i) 
$$
\left| \int_{\partial M} \, B_g \, \varphi_{a,\l} \, \varphi_{a,\l} d\s_g \quad - \ov{c}^{\frac{2(n - 1)}{n - 2}} \, \int_{\R^{n - 1}}\quad \frac{dx}{(1 + |x|^2)^{n - 1}} \right| \quad \leq  \, \frac{C}{\l^{n - 2}}\quad  \mbox{for } a \in \partial M 
$$
\item(ii) 
$$
\left| \int_{\partial M} \,  \varphi_{a,\l}^{\frac{2(n - 1)}{n - 2}}  d\s_g \quad - \ov{c}^{\frac{2(n - 1)}{n - 2}} \, \int_{\R^{n - 1}}\quad \frac{dx}{(1 + |x|^2)^{n - 1}} \right| \quad \leq  \, \frac{C}{\l^{n - 2}}\quad  \mbox{for } a \in \partial M 
$$
\item(iii)
$$
 \int_{\partial M} \,  \varphi_{a_1,\l}^{\frac{n}{n - 2}} \,  \varphi_{a_2,\l} d\s_g \quad \geq  \, \frac{C}{\l^{n - 2}} \quad  \mbox{for } a_1, a_2  \in \partial M
$$
\item(vi)
$$
\int_{\partial M} \, B_g \, \varphi_{a,\l} \, \varphi_{a,\l} d\s_g \quad \leq (1 + \theta) \, \int_{\partial M} \,  \varphi_{a_1,\l}^{\frac{n}{n - 2}} \,  \varphi_{a_2,\l} d\s_g
$$ 

\end{description}

\end{lem}

\begin{pf}

Proof of (i)

From the definition of $ \varphi_{a,\l}$, we derive:
$$
\int_{\partial M} \, B_g \, \varphi_{a,\l} \, \varphi_{a,\l} d\s_g = \\
\int_{\partial M} \,(\o_a \, \d_{a,\l} u_a)^{\frac{n}{n - 2}} \, \varphi_{a,\l} d\s_g
$$
Using Lemma \ref{l:error}, we deduce:
\begin{eqnarray}
\int_{\partial M} \, B_g \, \varphi_{a,\l} \, \varphi_{a,\l} d\s_g 
& =& \, \int_{\partial M \cap B_{2 \rho}} \,(\o_a \, \d_{a,\l} u_a)^{\frac{2(n - 1)}{n - 2}} \,  d\s_g 
 + \, O(\frac{1}{\l^{n - 2}}) \int_{\partial M \cap B_{2 \rho}} \,(\o_a \, \d_{a,\l} u_a)^{\frac{n}{n - 2}} \,  d\s_g \nonumber  \\
&= & \, \int_{\partial M \cap B_{ \rho}} \,( \d_{a,\l} )^{\frac{2(n - 1)}{n - 2}} \,  dv
    +   \, O(\frac{1}{\l^{n - 2}}) \int_{\partial M \cap B_{2 \rho}} \,(\o_a \, \d_{a,\l} u_a)^{\frac{n}{n - 2}} \,  dv_{g_0} + O(\l^{n - 1})\nonumber \\
& = & \ov{c}^{\frac{2(n - 1)}{n -2}} \quad \int_{\R^{n -1}} \, \frac{dx}{(1 + |x|^2)^{n - 1}} \, + O(\frac{1}{\l^{n - 2}})
\end{eqnarray}.

The proof of (ii) is essentially reduced, up to minor differences to the same computations involved in the proof of (ii).

Proof of (iii)

From Lemma \ref{l:lb} we deduce
$$
\int_{\partial M} \, \varphi_{a_1,\l}^{\frac{n }{n - 2}}\, \varphi_{a_2, \l} d\s_g \\ \geq \frac{1}{C\l^{n - 2}}\int_{\partial M} \, \varphi_{a_1,\l}^{\frac{n }{n - 2}}d\s_g
$$

Then from Lemma \ref{l:error} and Lemma \ref{l:lb} we derive
\begin{eqnarray*}
\int_{\partial M} \, \varphi_{a_1,\l}^{\frac{n }{n - 2}}\, \varphi_{a_2, \l} 
 d\s_g \, \geq \,  \frac{1}{C\l^{n - 2}}\int_{ B_{\rho}(a) \cap \partial M} \, \hat{\d}_{a_1,\l}^{\frac{n }{n - 2}}d\s_g \geq  \frac{1}{C\l^{n - 2}}\int_{B_{\rho}(a)\cap \partial M} \, \d_{a_1,\l}^{\frac{n }{n - 2}}d\s_{g_0} = O(\frac{1}{\l^{n - 2}})
\end{eqnarray*}

Proof of (vi)
\begin{eqnarray*}
\int_{\partial M} \, B_g\varphi_{a_1,\l}\, \varphi_{a_2, \l} 
 dv_g \,& = & \int_{\partial M} \, \hat{\d}_{a_1,\l}^{\frac{n }{n - 2}}\, \varphi_{a_2, \l} dv_g \,  \\
& = &\int_{\partial M \cap B_{\rho}} \, \varphi_{a_1,\l}^{\frac{n }{n - 2}} \var_{a_2,\l} \,+ O(\frac{1}{\l^{\frac{n}{2}}}) \int_{\partial M \setminus B_{\rho}}\varphi_{a_2, \l} \\
& = & \int_{\partial M} \, \varphi_{a_1,\l}^{\frac{n }{n - 2}}\, \varphi_{a_2, \l} d\s_g \,   \, + O(\frac{1}{\l^{\frac{n -2}{2}}}) \int_{\partial M} \, \hat{\d}_{a_1,\l}^{\frac{2 }{n - 2}}\, \varphi_{a_2, \l} d\s_g \, + O(\frac{1}{\l^{n - 1}}) \, \\
& = & \int_{\partial M} \, \varphi_{a_1,\l}^{\frac{n }{n - 2}}\, \varphi_{a_2, \l} dv_g \, \varphi_{a_2, \l}  \, + O(\l^{\frac{n -2}{2}}) \int_{\partial M} \, \d_{a_1,\l}^{\frac{2 }{n - 2}}\, \varphi_{a_2, \l} d\s_{\ov{g}_{0}}  \, \\
& = & \int_{\partial M} \, \varphi_{a_1,\l}^{\frac{n }{n - 2}}\, \varphi_{a_2, \l} dv_g \,   \, + O(\frac{1}{\l^{\frac{n -2}{2}}}) \int_{\partial M} \, \d_{a_1,\l}^{\frac{2 }{n - 2}}\, \d_{a_2, \l} d\s_{\ov{g}_{0}} \, + O(\frac{1}{\l^{n - 1}}) \,
\end{eqnarray*}

Now from  Lemma \ref{a1} in the Appendix we deduce :
\begin{eqnarray*}
O(\frac{1}{\l^{\frac{n -2}{2}}}) \int_{\partial M \cap B_{\rho}(a_1)} \, \d_{a_1,\l}^{\frac{n }{n - 2}}\, \d_{a_2, \l} d\s_{\ov{g}_{0}} \, & = & o \left( \int_{\partial M} \, \d_{a_1,\l}^{\frac{n }{n - 2}}\, \d_{a_2, \l} d\s_{\ov{g}_{0}} \right) 
\end{eqnarray*}

Therefore using (iii) we have
$$
\int_{\partial M} \, B_g\varphi_{a_1,\l}\, \varphi_{a_2, \l} 
 dv_g \, =\int_{\partial M} \, \varphi_{a_1,\l}^{\frac{n}{n - 2}}\, \varphi_{a_2, \l}dv_g  (1 + o(1))
$$
The proof of (vi) and the proof of Lemma \ref{l:fe} are thereby completed.
\end{pf}

\section{ Some standard facts }

We recall that solutions of Problem $(P)$ arises , up to a constant, as critical points of the   functional $ J$ is defined by

$$
J(u) = \left( \int_{M} \, -L_gu \, u \, dv_g \, + \, \int_{\partial M} \, B_g u \, u \, d\s_g \right)^{\frac{n - 1}{n - 2}} \left({\int_{\partial M} \, u^{\frac{2(n - 1)}{n - 2}}}\right)^{ - 1}
$$
 where $u$ belongs to  $\Sigma^+$ defined as follows:
$$
\Sig^+ = \{ u  \in H^1(M), u \geq 0, \barre{u} \, = 1 \}
$$
Let us observe that $ \Sig^+$ is invariant by the flow of $- \partial J $.

The functional $J$ is known to not satisfy Palais Smale condition(PS for short) , which leads to the failure of classical existence mecanism. In order to describe this failure we need some notation.

For $\e > 0$ and $p \geq 1$ , let 
$$
V(p,\e)   =    \left\{ \begin{array}{cc}
   u \in \Sig^+  \, \mbox{such that } \exists \,  (a_1,\cdots,a_p) \in (\partial M)^p \, \mbox{and } \exists \, (\l_1,\cdots,\l_p) \in (\R^*_+)^p \, \mbox{such that }   & \hskip -.4truecm\\
\salt
\barre{u - \frac{\sum_{i = 1}^p \var_{a_i,\l_i}}{ ||\sum_{i = 1}^p \var_{a_i,\l_i}  ||}} < \e , \, \mbox{with } \, \l_i \geq \frac{1}{\e} \, \mbox{and } \e_{ij} < \e \, &  \hskip -.4truecm
\end{array}
\right\}.
$$
where $ \e_{ij} \, = \, \left( \frac{1}{\frac{\l_i}{\l_j} + \frac{\l_j}{\l_i} + \l_i \l_j \, d(a_i , a_j)^2}  \right)^{\frac{n - 2}{2}}$
and $d$ denotes the geosedic distance.

If $u$ is a function in $V(p,\e) $ , one can find an optimal representation, arguing as in Proposition 7  of \cite{bc} , namely  we have:

\begin{lem}\label{l:rep}
For every $p \geq 1$, there exists $ \e > 0$ such that $ \forall \, u \in V(p,\e)$ the minimization problem
$$
\inf_{\a_i,b_i,\mu_i} \, \barre{u - \sum_{i = 1}^p \a_i \var_{b_i,\mu_i}}
$$
has a unique solution, up to permutation on the set of indices $ \{1,\cdots,p\}$
\end{lem}
At this point we introduce the following notations:

 Let $S$ be defined as $ S = \left( \int_{\R^{n -1}}\, \frac{1}{(1 + |x|^2)^{n - 1}}  \right)^{\frac{1}{n - 2}}$,   $ b_p = p^{\frac{1}{n - 2}} S $ 
  and 
$$
 W_p = \{ u \in \Sig^+ , \mbox{such that } J(u) < b_{p + 1} \}
$$
We are ready now to state the characterization of the Palais Smale sequences failing the P.S condition.

\begin{pro}\label{p:ps}
Under the assumption that (P) has no solution, 
let $ u_k \,  \subset \Sig^+$ be a sequence satisfying $ J(u_k) \to c $ , a positive number and $ \partial J (u_k) \to 0$ . There exist an integer $ p \geq 1$ and a sequence $ (\e_k)_k$ such that $ u_k \in V(p,\e)$.
Conversely, let $ p \in \N^+$, let $ \e_k$ be a positive sequence with $ \lim_{k \to + \infty} \e_k = 0$ and let $ u_k \in V(p,\e)$ then $\partial J(u_k) \to 0$ and $ J(u_k) \to b_p$.
\end{pro}

\begin{pf}
The proof of this Proposition is by now standard, taking into account the uniqueness result of Li-Zhu \cite{LZ1}, see also \cite{E3} and using the  Liouville Theorem to rule out the possibility of interior blow up.
\end{pf}
we have also the following  local deformation Lemma , similar to Lemma 17 in \cite{bb} :

\begin{lem}\label{l:def1}
Under the asumption that (P) has no solution, for $ \e > 0$, the pair $ (W_p,W_{p - 1})$ retracts by deformation onto the pair $ (W_{p - 1} \cup A_p, W_{p - 1}) $ where $ A_p \subset V(p,\e) $.

\end{lem}

\section{Expansion of the functional near its potential critical point at infinity}

This section is devoted to an asymptotic expansion of $J$ in the neighborhood of its potential critical points at infinity, that is in some $V(p,\e)$. This expansion displays the fact that when the number of the bubbles are large enough, their interaction increases to force their energy to be under their critical level.Such a fact is a key point in the topological argument.

\begin{lem}\label{l:be}
There holds:

\begin{description}
\item(i) For every $p \in \N^*$ and every $ \e_1 > 0$, there exists $ \l_p = \l(p, \e_1)$ such that for any $ (\a_1, \cdots,\a_p)$ satisfying  $ \a_i \geq 0, \sum_{i = 1}^p \, \a_i \, = 1$ , for any $ (x_1, \cdots,x_p) \in \,  (\partial M)^p$ for any $ \l \geq \l_p$, we have:
$$
 \mbox{If } \sum_{i \not= j} \, \int_{\partial M} \, \varphi_{a_i,\l}^{\frac{n}{n - 2}} \, \varphi_{a_j,\l}  \geq \e_1 \quad \mbox{then } \quad J(\sum_{i = 1}^p \, \a_i \varphi_{\a_i, \l}) \, \leq \, p^{\frac{1}{n - 2}}S .
$$
\item(ii) There exist $ 0 < \theta_0 < 1$, $C > 0$, and $ \ov{\e_1 } \geq \e_1$, such that for any $p \in \N^*$, for any $ (\a_1, \cdots, \a_p)$ satisfying $ \a_i \geq 0,\quad  \sum_{i = 1}^p \a_i = 1 , \quad \frac{\a_i}{\a_j} \geq \theta_0$ , for any $ (a_1, \cdots, a_p) \in (\partial M)^p$, for any $ \l \geq 1$, the following inequality holds
$$
J(\sum_{i = 1}^p \, \a_i \varphi_{\a_i, \l}) \, \leq \, p^{\frac{1}{n - 2}} \,
 S\left( 1 +  O(\frac{1}{\l^{n - 2}}) \, + \quad \frac{(p + 1)C}{\l^{n - 2}}\right)
$$
If 
$$
 \sum_{i \not= j} \, \int_{\partial M} \, \varphi_{a_i,\l}^{\frac{n}{n - 2}} \, \varphi_{a_j,\l} \leq \ov{\e}_1
$$

\item(iii) If we drop in (ii) the condition $\frac{\a_i}{\a_j} \geq \theta_0 $, then the following weaker inequality still holds:
$$
J(\sum_{i = 1}^p \, \a_i \varphi_{\a_i, \l}) \, \leq \, p^{\frac{1}{n - 2}} S( 1 + \, \frac{1}{C\l^{n - 2}} + \quad \frac{1}{C}  \sum_{i \not= j} \, \int_{\partial M} \, \varphi_{a_i,\l}^{\frac{n}{n - 2}} \, \varphi_{a_j,\l}  
$$
\end{description}
\end{lem}

\begin{pf}

Let 
$$
\begin{cases}
b_i  = \, \frac{\a_i \hat{\d}_{a_i,\l}}{\sum_{j = 1}^p \, \a_j \, \hat{\d}_{a_j,\l}} \quad \mbox{if } \hat{\d}_{a_i,\l} \not= 0  \, \\
 b_i  = \quad  1  \quad  \mbox{if }  \hat{\d}_{a_i,\l} = 0

\end{cases}
$$
then Lemma B2  of \cite{bc}, which extends to our functional,  implies that

$$
J(\sum_{i = 1}^p \, \varphi_{a_i,\l}) \, \leq \,\left( \frac{\int_{\partial M} \, (\sum_{i = 1}^p\, \hat{\d}_{a_i,\l})^{\frac{2(n - 1)}{n - 2}}}{\int_{\partial M}(\sum_{i = 1}^p \, \varphi_{a_i,\l})^{\frac{2(n - 1)}{n - 2}}}  \right)^{\frac{1}{2}}  \quad \left( \sum_{i = 1}^p \, b_i \, \int_{\partial M} \, \hat{\d}_{a_i,\l}^{\frac{2(n - 1)}{n - 2}}\right)^{\frac{1}{n - 2}}
$$
Using Lemma \ref{l:error}, and $\int_{\partial M} \, \hat{\d}_{i}^{\frac{n}{n - 2}}\, dv_g \quad = \, O(\frac{1}{\l^{\frac{n - 2}{2}}}) $, we derive easily

$$
J(\sum_{i = 1}^p \, \varphi_{a_i,\l})  \leq \, (1 + O(\frac{1}{\l^{n - 2}})) \quad \left( \sum_{i = 1}^p \, b_i \, \int_{\partial M} \, \hat{\d}_{a_i,\l}^{\frac{2(n - 1)}{n - 2}}\right)^{\frac{1}{n - 2}}
$$
Thus we  obtain
$$
J(\sum_{i = 1}^p \, \varphi_{a_i,\l})  \leq \, (1 + O(\frac{1}{\l^{n - 2}})) \quad \left(   (p - 1)S^{n - 2} \, + \, \int_{\partial M} \, \frac{\a_1\, \hat{\d}_{a_1,\l}}{\a_1 \, \hat{\d}_{a_1,\l} \, + \, \a_2\hat{\d}_{a_2,\l}} \hat{\d}_{a_1,\l}^{\frac{2(n - 1)}{n - 2}}\right)^{\frac{1}{n - 2}}
$$
We may assume without loss of generality that 

\begin{description}
\item(i) $ \frac{\a_1}{\a_2} \leq 1$
\item(ii)
$$
  \int_{\partial M} ( \varphi_{a_1,\l}^{\frac{n}{n - 2}} \varphi_{a_2,\l} \, + \,  \varphi_{a_1,\l}  \varphi_{a_2,\l}^{\frac{n}{n - 2}}) \,  =  \sup_{i \not= j} \int_{\partial M}  \varphi_{a_i,\l}^{\frac{n}{n - 2}} \varphi_{a_j,\l} \, + \,  \varphi_{a_i,\l}  \varphi_{a_j,\l}^{\frac{n}{n - 2}} \,  
$$
\end{description}

We then have using that $ \int_{\partial M \setminus B_{\frac{\rho}{2}}(a_1)} \hat{\d}_{a_1,\l}^{\frac{2(n - 1)}{n - 2}} dv_g  \, = O(\frac{1}{\l^{n - 2}})$

$$
J(\sum_{i = 1}^p \, \varphi_{a_i,\l})  \leq \, (1 + O(\frac{1}{\l^{n - 2}})) \quad \left(   (p - 1)S^{n - 2} \, + \, \int_{\partial M} \, \frac{\a_1\, \hat{\d}_{a_1,\l}}{\a_1 \, \hat{\d}_{a_1,\l} \, + \, \a_2\hat{\d}_{a_2,\l}} \d_{a_1,\l}^{\frac{2(n - 1)}{n - 2}}\right)^{\frac{1}{n - 2}}
$$
The continuity of $u_y$ with respect to $y$ implies the existence of $ \eta > 0$ such that
$$
\frac{1}{2} \, \leq \,  \frac{u_{a_1}}{u_{a_2}} \leq 2 \quad \mbox{if } \quad x_2 \in B_{\eta}(a_1)
$$
Thus if $ d(a_2,a_2) \, \leq \, \eta$, we have

\begin{eqnarray}\label{e:m1}
J(\sum_{i = 1}^p \, \varphi_{a_i,\l})  
 & \leq  \, & (1 + O(\frac{1}{\l^{n - 2}})) \quad \left(   (p - 1)S^{n - 2} \, + \, \int_{\partial M} \, \frac{\a_1\, \hat{\d}_{a_1,\l}}{\a_1 \, \hat{\d}_{a_1,\l} \, + \, \a_2\hat{\d}_{a_2,\l}} \hat{\d}_{a_1,\l}^{\frac{2(n - 1)}{n - 2}}\right)^{\frac{1}{n - 2}} \, \nonumber  \\
& \leq & (1 + O(\frac{1}{\l^{n - 2}})) \quad \left(   (p - 1)S^{n - 2} \, - \, \int_{\partial M \cap B_{a_1}(\rho)} \, \frac{\a_1\, \d_{a_1,\l}}{\a_1 \, \hat{\d}_{a_1,\l} \, + \,  \d_{a_2,\l}} \hat{\d}_{a_1,\l}^{\frac{2(n - 1)}{n - 2}}\right)^{\frac{1}{n - 2}} \, 
\end{eqnarray}.

At this point we state the following Claim  which proof is postponed until the end of this section.

{\bf Claim} : There exists $\e_0$ small enough, and a positive constant  $\ov{C}'$ such that
$$
\int_{\partial M} \frac{\d_{a_2,\l}}{ 2 \d_{a_2,\l} \, + \, \d_{a_2,\l}} \, \d_{a_1,\l}^{\frac{2(n - 1)}{n - 2}}   \, dv_{\ov{g}_{0}} \geq \ov{C}' \e_0
\, \int_{\partial M} \varphi_{a_1,\l}^{\frac{n}{n - 2}}\varphi_{a_2,\l} \, + \, \varphi_{a_1,\l} \, \varphi_{a_2,\l}^{\frac{n}{n - 2}} dv_{\ov{g}_{0}}
$$
From another part , by assumption, we know that

$$
\sum_{i \not= j} \, \int_{\partial M} \left( \varphi_{a_i,\l}^{\frac{n}{n - 2}} \varphi_{a_j,\l} dv_{\ov{g}} \right) \, \leq \, p^2 \, \int_{\partial M} \left( \varphi_{a_1,\l}^{\frac{n}{n - 2}} \varphi_{a_2,\l}  +   \varphi_{a_2,\l}^{\frac{n}{n - 2}} \varphi_{a_1,\l} \right) dv_{\ov{g}}
$$
Thus we derive from the above claim and (\ref{e:m1}) , that: 

\begin{eqnarray}\label{eq:m2}
J(\sum_{i = 1}^{p} \a_i \varphi_{a_i,\l}) & \leq & (1 + O(\frac{1}{\l^{n - 2}})) \left( pS^{n - 2} - \, \frac{\ov{C} \e_0}{p^2} \sum_{i \not= j} \int_{\partial M} \varphi_{a_i,\l}^{\frac{n}{n - 2}} \, \varphi_{a_j,\l} dv_{\ov{g}}
\right)^{\frac{1}{n - 2}} \, \\
& \leq &  p^{\frac{1}{n - 2}} S (1 + O(\frac{1}{\l^{n - 2}}) - \frac{\ov{C_1} \e_0 \e_1}{p^3})
\end{eqnarray}
clearly implies (i), which is therefore proven if $ d(a_1,a_2) \leq \eta $.
Now we rule out the case where $ d(a_1,a_2) \geq \eta$ as follows:

If $ d(a_1,a_2) \geq \eta$ then
$$
 \int_{\partial M} \left(\varphi_{a_1,\l}^{\frac{n}{n - 2}} \, \varphi_{a_2,\l} \, + \,\varphi_{a_2,\l}^{\frac{n}{n - 2}} \, \varphi_{a_1,\l} \right) \, = \, O(\frac{1}{\l^{n - 2}})
$$ Therefore
$$
\sum_{i \not= j}\varphi_{a_i,\l}^{\frac{n}{n - 2}} \, \varphi_{a_j,\l} \, \leq p^2 .O(\frac{1}{\l^{n - 2}})  \, = \, O(\frac{1}{\l^{n - 2}}) 
$$
Taking $ \l$ very large we have $\sum_{i \not= j}\varphi_{a_i,\l}^{\frac{n}{n - 2}} \, \varphi_{a_j,\l} \, < \e_1 $ therefore the proof of (i) is reduced to the case $ d(a_1,a_2) \leq \eta$.
The proof of (i ) is thereby established.

The proofs of (ii) and (iii) will be completed toghether since they rest on the same expansion of  the functional $J$.

Using Lemma \ref{l:ineq} we derive the following inequality
\begin{equation}\label{eq:e1}
 J(\sum_{i - 1}^p \a_i \, \varphi_{a_i,\l})\, 
 \leq  \frac{\left(\sum_{i = 1}^p \a_i^2 \, \int_{\partial M} \hat{\d}^{\frac{n}{n - 2}} \varphi_{a_i,\l} \, + \, \sum_{i \not= j}^p \a_i \a_j \int_{\partial M} \hat{\d}_{a_i,\l}^{\frac{n}{n - 2}} \varphi_{a_j,\l} \right)^{\frac{n - 1}{n - 2}}}{ \sum_{i = 2}^p \a_i^{\frac{2(n - 1)}{n - 2}} \, \int_{\partial M} \varphi_{a_i,\l}^{\frac{2(n - 1)}{n - 2}} \, + \, \sum_{i \not= j}^p \frac{\g(n - 1)}{n - 2} \a_i^{\frac{n}{n - 2}} \a_j \, \int_{\partial M} \varphi_{a_i,\l}^{\frac{n}{n - 2}} \varphi_{a_j,\l} }
\end{equation}

Then (\ref{eq:e1}) implies using Lemma \ref{l:error} and Lemma \ref{l:fe}

\begin{eqnarray}\label{eq:e2}
 & {} & J(\sum_{i - 1}^p \a_i \, \varphi_{a_i,\l})\, \\
&  \leq  & \frac{\left(\sum_{i = 1}^p \a_i^2 \, S^{n - 2} \, +  \, + \,( 1 + \theta) \sum_{i \not= j}^p \a_i \a_j \int_{\partial M} \varphi_{a_i,\l}^{\frac{n}{n - 2}} \varphi_{a_j,\l} \ + \ \sum_{i = 1}^p \a_i^2 \frac{\ov{C}}{\l^{n - 2}} \right)^{\frac{n - 1}{n - 2}}}{ \sum_{i = 2}^p \a_i^{\frac{2(n - 1)}{n - 2}}  S^{n - 2} \, + \, \sum_{i \not= j}^p \frac{\g(n - 1)}{n - 2} \a_i^{\frac{n}{n - 2}} \a_j \, \int_{\partial M} \varphi_{a_i,\l}^{\frac{n}{n - 2}} \varphi_{a_j,\l}  \, - \, \sum_{i = 1}^p \a_i^{\frac{2(n - 1)}{n - 2}} \frac{\ov{C}}{\l^{n - 2}}} \nonumber
\end{eqnarray}

where $\ov{C}$ and $ \ov{C}'$ are positive constants independant of $ (\a_1,\cdots,\a_p)$, $ \l$ and $p$.

Let us assume that 
$$
\sum_{i \not= j} \int_{\partial M} \varphi_{a_i,\l}^{\frac{n}{n - 2}} \varphi_{a_j,\l} \, < \, \e_1.
$$
If $ \e_1$ is chosen small enough, then for $\l$ large enough , we have

\begin{eqnarray}\label{eq:e3}
 & {} & J(\sum_{i - 1}^p \a_i \, \varphi_{a_i,\l})\, \nonumber \\
&  \leq  & \frac{\left(\sum_{i = 1}^p \a_i^2 \, S^{n - 2} \, +  \, + \,( 1 + \theta) \sum_{i \not= j}^p \a_i \a_j \int_{\partial M} \varphi_{a_i,\l}^{\frac{n}{n - 2}} \varphi_{a_j,\l} \ + \ \sum_{i = 1}^p \a_i^2 \frac{\ov{C}}{\l^{n - 2}} \right)^{\frac{n - 1}{n - 2}}}{ \sum_{i = 2}^p \a_i^{\frac{2(n - 1)}{n - 2}}  S^{n - 2} \, + \, \sum_{i \not= j}^p \frac{\g(n - 1)}{n - 2} \a_i^{\frac{n}{n - 2}} \a_j \, \int_{\partial M} \varphi_{a_i,\l}^{\frac{n}{n - 2}} \varphi_{a_j,\l}  \, - \, \sum_{i = 1}^p \a_i^{\frac{2(n - 1)}{n - 2}} \frac{\ov{C}}{\l^{n - 2}}}
\end{eqnarray}.
Making an expansion we have
\begin{eqnarray}\label{eq:e2a}
 & {} & J(\sum_{i - 1}^p \a_i \, \varphi_{a_i,\l})\, \\
&  \leq  & \frac{\left(\sum_{i = 1}^p \a_i^2 \, S^{n - 2} \, +  \, + \,( 1 + \theta) \sum_{i \not= j}^p \a_i \a_j \int_{\partial M} \varphi_{a_i,\l}^{\frac{n}{n - 2}} \varphi_{a_j,\l} \ + \ \sum_{i = 1}^p \a_i^2 \frac{\ov{C}}{\l^{n - 2}} \right)^{\frac{n - 1}{n - 2}}}{ \sum_{i = 2}^p \a_i^{\frac{2(n - 1)}{n - 2}}  S^{n - 2} \, + \, \sum_{i \not= j}^p \frac{\g(n - 1)}{n - 2} \a_i^{\frac{n}{n - 2}} \a_j \, \int_{\partial M} \varphi_{a_i,\l}^{\frac{n}{n - 2}} \varphi_{a_j,\l}  \, - \, \sum_{i = 1}^p \a_i^{\frac{2(n - 1)}{n - 2}} \frac{\ov{C}}{\l^{n - 2}}}
\end{eqnarray}
Finally we obtain:

 \begin{eqnarray*}\label{eq:e4} 
 J(\sum_{i - 1}^p \a_i \, \varphi_{a_i,\l})\, 
&  \leq  & \frac{\left(\sum_{i = 1}^p \a_i^2 \right)^{\frac{n - 1}{n - 2}} }{\sum_{i = 1}^p \a_i^{\frac{2(n - 1)}{n - 2}}} S 
\left( 1 +  \frac{n - 1}{n - 2} \, \sum_{i \not= j} \frac{\a_i \a_j}{\sum_i \a_i^2 S^{n - 2}}  \, \int_{\partial M} \var_{a_i,\l}^{\frac{n}{n - 2}} \var_{a_j,\l} \,
 \right.\\
 & { } & \left. \, - \, \frac{n - 1}{n - 2} \, \sum_{i \not= j} \frac{\a_i \a_j}{\sum_i \a_i^2 S^{n - 2}}  \, \int_{\partial M} \var_{a_i,\l}^{\frac{n}{n - 2}} \var_{a_j,\l} + \frac{\ov{C}''}{\l^{n - 2}}  \right)
\end{eqnarray*}
Let us observe that  (iii) follows from (\ref{eq:e4}), so it remains only to prove (ii).

Proof of (ii)

Let us now assume that there exists $\th_0\, , \, 0 < \th_0 < 1 $ and $ \frac{\a_i}{\a_j} \geq \th_0$ for any $(i,j)$, then
$$
\frac{\a_i \a_i}{\sum_{r = 1}^p \a_r^2} \leq \frac{1}{p \th_0^2} \quad  \mbox{and } \quad \frac{\a_i^{\frac{n}{n - 2}} \a_j}{\sum_{r = 1}^p \a_r^{\frac{2(n - 1)}{n - 2}}} \, \geq \frac{\th_0^{\frac{2(n - 1)}{n - 2}}}{p}
$$
Now we choose $ 0 < \th_0 < 1$ and $ \th > 0$ such that $ (1 + \th)\frac{1}{\th_0^2} \, - \, \g \th_0^{\frac{2(n - 1)}{n - 2}} < 0$.

Let $ \d = S^{n - 2} \left( \frac{\g(n - 1)}{n - 2} \th_0^{\frac{2(n - 1)}{n - 2}} \, - \,\frac{n - 1}{n - 2} \frac{1 + \th}{\th_0^2}  \right)$

We then derive 
\begin{equation}\label{eq:ef} 
J(\sum_{i - 1}^p \a_i \, \varphi_{a_i,\l}) \, \leq \,\frac{\left(\sum_{i = 1}^p \a_i^2 \right)^{\frac{n - 1}{n - 2}} }{\sum_{i = 1}^p \a_i^{\frac{2(n - 1)}{n - 2}}} S 
\left( 1 \, - \, \frac{\d}{p} \, \sum_{i \not= j} \frac{\a_i \a_j}{\sum_i \a_i^2 S^{n - 2}}  \, \int_{\partial M} \var_{a_i,\l}^{\frac{n}{n - 2}} \var_{a_j,\l} + \frac{\ov{C}''}{\l^{n - 2}}  \right)
\end{equation}

Then using (iii) of Lemma \ref{l:fe} , we derive (ii) from (\ref{eq:ef}).
The proof of Lemma \ref{l:be} is thereby complete.\end{pf}

{\bf Proof of the claim }

Let $\e_0 > 0$ be given and let 
$$
E_0 = \{x \in \partial M \quad \mbox{such that }\d_{a_1,\l}(x) \geq \e_0 \left(  2 \d_{a_1,\l} \, + \, \d_{a_2,\l} \right)(x)   \}
$$
We have:
\begin{eqnarray}\label{eq:cl1}
& { } & \int_{\partial M} \frac{\d_{a_2,\l}}{2 \d_{a_1,\l} + \d_{a_2,\l}} \, \\
& \geq & \e_0 \, \int_{E_0} \d_{a_1,\l}^{\frac{n}{n - 2}} \, \d_{a_2,\l} \,  dv_{\ov{g}_{0}} \,\nonumber \\
& \geq  & \e_0 \, \left( \int_{\partial M} \d_{a_1,\l}^{\frac{n}{n - 2}} \, \d_{a_2,\l} \, dv_{\ov{g}_{0}}  \, - \, (\frac{\e_0}{1 - 2\e_0})^{\frac{1}{n - 2}} \int_{\partial M} \d_{a_1,\l}^{\frac{n - 1}{n - 2}} \, \d_{a_2,\l}^{\frac{n - 1}{n - 2}}\right) \nonumber
\end{eqnarray}
Let us obseve that 
$$
\int_{\partial M} \d_{a_1,\l}^{\frac{n}{n - 2}} \, \d_{a_2,\l} = \, \int_{\partial M} \d_{a_2,\l}^{\frac{n}{n - 2}} \, \d_{a_1,\l}
$$
and
$$
\int_{\partial M} \d_{a_1,\l}^{\frac{n - 1}{n - 2}} \, \d_{a_2,\l}^{\frac{n - 1}{n - 2}} \leq \frac{1}{2} \left( \int_{\partial M} \d_{a_1,\l}^{\frac{n}{n - 2}} \d_{a_2,\l} \, + \, \d_{a_1,\l} \d_{a_2,\l}^{\frac{n}{n - 2}}  \right)
$$
Then (\ref{eq:cl1}) becomes

\begin{equation}\label{eq:cl2}
 \int_{\partial M} \frac{\d_{a_2,\l}}{2 \d_{a_1,\l} + \d_{a_2,\l}} \, \geq \frac{\e_0}{2} \, (1 - (\frac{\e_0}{1 - \e_0})^{\frac{1}{n - 2}}) \,  \int_{\partial M} \d_{a_1,\l}^{\frac{n}{n - 2}} \d_{a_2,\l} \, + \, \d_{a_1,\l} \d_{a_2,\l}^{\frac{n}{n - 2}}
\end{equation}
Thus for $\e_0$ small enough , we have:
\begin{eqnarray}\label{cl3}
 \int_{\partial M} \frac{\d_{a_2,\l}}{2 \d_{a_1,\l} + \d_{a_2,\l}} \, & \geq  & C \e_0  \int_{\partial M} \d_{a_1,\l}^{\frac{n}{n - 2}} \d_{a_2,\l} \, + \, \d_{a_1,\l} \d_{a_2,\l}^{\frac{n}{n - 2}}  \nonumber, \\
& \geq &  C' \e_0 \, \int_{\partial M \cap B_{\rho}(a_1)} \hat{\d}_{a_1,\l}^{\frac{n}{n - 2}} \ov{\d}_{a_2,\l} \, + \, \ov{\d}_{a_1,\l} \hat{\d}_{a_2,\l}^{\frac{n}{n - 2}}  \nonumber \, \\
& \geq & \,  C'' \e_0 \,   \int_{\partial M \cap B_{\rho}(a_1)} \var_{a_1,\l}^{\frac{n}{n - 2}} \var_{a_2,\l} \, + \, \var_{a_1,\l} \var_{a_2,\l}^{\frac{n}{n - 2}} \, \nonumber \\
& \geq & \,  C'' \e_0 \,   \int_{\partial M } \var_{a_1,\l}^{\frac{n}{n - 2}} \var_{a_2,\l} \, + \, \var_{a_1,\l} \var_{a_2,\l}^{\frac{n}{n - 2}} .
\end{eqnarray}
 Hence our claim is proved.

Lemma \ref{l:be} implies the following proposition:

\begin{pro}\label{p:bp}
There exists an integer $ p_0$ and a positive real number $ \l_0 > 0$ such that  for any $ (\a_1, \cdots, \a_p)$ satisfying $ \a_i \geq 0, \, \sum_{i = 1}^p
\, \a_i = 1$, for any $ (a_1,\cdots,a_p) \in \partial M$, for any $ \l \geq \l_0$, we have
$$
J(\sum_{i = 1}^p \, \a_i \varphi_{\a_i, \l}) \, \leq \, p_0^{\frac{1}{n - 2}} S
$$

\end{pro}

\begin{pf}
The proof of Proposition \ref{p:bp} follows from (i), (ii) and  (iii) of Lemma \ref{l:be}.
We first choose $ 0 < \e_1 < \ov{\e_1} $ and $ \l_0$ so that :
\begin{equation}\label{eq:n1}
\frac{(\sum_{i = 1}^p \a_i^2)^{\frac{n - 1}{n - 2}}}{\sum_{i = 1}^p \a_i^{\frac{2(n - 1)}{n - 2}} } \left( 1 \, + \, O(\frac{1}{\l^{n - 2}}) \, + \, \frac{\e_1}{C}   \right) \, < p^{\frac{1}{n - 2}}
\end{equation} .

Considering $(\a_1,\cdots,\a_p) $ , $(a_1,\cdots,\a_p) $ and $ \l \geq \l_0$ , we study various cases:

\begin{description}
\item[1st case]:
There exists $ (i_0,j_0)$ such that $\frac{\a_{i_0}}{\a_{j_0}} \, \leq \, \th_0 $ , then taking, $\l \geq \l_{p_0} = \sup(\l(p_0,\e),\l_0)$ where $ \l(p_0,\e_1)$ is given by (i) of Lemma \ref{l:be} , we derive :
$$
J(\sum_{i = 1}^p \a_i \, \var_{a_i,\l}) \, \leq \, p^{\frac{1}{n - 2}} S \quad \mbox{if } \, \sum_{i \not= j}^p \int_{\partial M} \var_{a_i,\l} \, \var_{a_j,\l} \, d\s_g \geq \e_1
$$
If on the contrary $ \sum_{i \not= j}^p \int_{\partial M} \var_{a_i,\l} \, \var_{a_j,\l} \, d\s_g \leq \e_1 $ we apply (iii). Since we have choose $ \e_1$ such that :

\begin{equation}\label{eq:nn}
\frac{(\sum_{i = 1}^p \a_i^2)^{\frac{n - 1}{n - 2}}}{\sum_{i = 1}^p \a_i^{\frac{2(n - 1)}{n - 2}} } \left( 1 \, + \, O(\frac{1}{\l^{n - 2}}) \, + \, \frac{\e_1}{C}   \right) \, < p^{\frac{1}{n - 2}}
\end{equation}
we derive that $J(\sum_{i = 1}^p \a_i \, \var_{a_i,\l}) \, \leq \, p^{\frac{1}{n - 2}} S $ and the proof of Proposition \ref{p:bp} is established in this case.
\item[2nd case]:
Let us assume that $\frac{\a_i}{\a_j} > \th_0 $ for any $ (i,j)$ , then either (i) or (ii) of Lemma \ref{l:be} holds. If (i) holds then Proposition \ref{p:bp} holds, so let assume that (ii) holds and then choose $ p_0$ such that : $ (p_0 + 1) c^2 > 1$ then 
$$
J(\sum_{i = 1}^p \a_i \, \var_{a_i,\l}) \, \leq \, p_0^{\frac{1}{n - 2}} S
$$
\end{description}

The proof of Proposition \ref{p:bp} is thereby complete.
\end{pf}

\section{Proof of Theorem \ref{t:mt}}

For the  proof of the Theorem \ref{t:mt} , we introduce the following notations:

For any $ p \geq 1$ and $ \l > 0$ , let 
$$
B_p = B_p(\partial M) = \{ \sum_{i = 1}^p \a_i \d_{a_{i}} \, , \a_i \geq 0, \, \sum_{i = 1}^p \a_i = 1 , \, a_i \in \partial M \, \}
$$
and $ B_0 = B_0(\partial M) = \emptyset $

Set also $ f_p(\l)$ to denote the map from $ B_p(\partial M)$ to $ \Sig^+$ defined by 
$$
f_p(\l) (\sum_{i = 1}^p \a_i \, \d_{a_{i}}) \, = \, \frac{\sum_{i = 1}^p \, \var_{a_{i},\l_i}}{\barre{\sum_{i = 1}^p \, \var_{a_{i},\l_i}}}
$$

Clearly we have $ B_{p - 1} \subset B_p$ and $ W_{p - 1} \subset W_p$.

Moreover $ f_p(\l)$ enjoys the following properties:

\begin{pro}\label{p:tt}
The function $ f_p(\l)$ has the following properties:
\begin{description}

\item(i) For any integer $ p \geq 1$ , there exists a real number $\l_p > 0 $ such that 
$$
f_p(\l) : \, B_p(\partial M) \to W_p \, \mbox{ for any } \l \geq \l_p
$$

\item(ii) 
There exists an integer $ p_0 > 1$ , such that for any integer $ p \geq p_0$ , and for any $ \l \geq \l_{p_{0}}$ , the map of pairs $ f_p(\l) : \, (B_p,B_{p - 1}) \, \to (W_p, W_{p - 1})$  satisfies $ (f_p)_*(\l) \equiv 0$
where 
$$
 (f_p(\l))_* : \,H_* (B_p,B_{p - 1}) \, \to H_*(W_p, W_{p - 1})
$$
and $H_* $is the $*$th homology group with $ \Z_{2} $ coefficients.
\end{description}
\end{pro}

\begin{pf}

(i) is a direct consequence of the inequalities (i), (ii) and  (iii) of Lemma \ref{l:be} , indeed  :
$$
J(f_p(\l)(\sum_{i = 1}^p \a_i \, \d_{a_{i}})) \, = \, J(\sum_{i = 1}^p \a_i \, \var_{a_{i},\l_i})
$$
(ii) follows from Propostion \ref{p:bp}.
\end{pf}

For the sequel we need the following notations:
Let $ \D_{p - 1} = \{ (\a_1, \cdots,\a_p), \a_i \geq 0, \sum_{i = 1}^p \a_i = 1\}$ and $ F_p = \{ (a_1,\cdots,a_p) \in (\partial M)^p \quad \mbox{such that } \exists \, i \not= j \quad \mbox{with } a_i = a_j \}$.
Let $ \s_p$ be the symmetric group of order p, which acts on $ F_p$, and let $T_p $ be a $ \s_p$- equivariant tubular neighborhood of $F_p$, in $(\partial M)^p $ (The existence of a such neighborhood is derived in the book of G. Bredon \cite{br})

From another part,  considering the topological pair $(B_p,B_{p - 1}) $ we observe that $(B_p \setminus B_{p - 1}) $ can be described as $ ((\partial M)^p)^* \times _{\s_p} (\D_p \setminus \partial \D_{p - 1}) $ where $(\partial M)^p)^* = \{ (a_1,\cdots,a_p) \in (\partial M)^p \quad \mbox{such that } a_i \not= a_j , \forall i \not= j \}  $
We notice that $ ((\partial M)^p)^*\times _{\s_p} (\D_p \setminus \partial \D_{p - 1})  $ is a noncompact manifold of dimension $ (n - 1)p + p - 1$.
Let for $ 0 < \theta < 1$,  $ \mathcal{M}_p = V_p \times _{\s_p} \D_{p - 1}^{\theta}$, where $ V_p = \ov{M^p \setminus T_p }$ and $\D_{p - 1}^{\theta} = \{ (\a_1,\cdots,\a_p) \in \D_{p - 1} \quad \mbox{such that } \frac{\a_i}{\a_j} \in [1 - \theta, 1 + \theta ], \forall i, \forall j \} $.
$ \mathcal{M}_p$ is a manifold which can be seen as a subst of $ B_p$, and the topological pair $(B_p,\mathcal{M}^c) $ retracts by deformation onto $ (B_p, B_{p - 1})$, we thus have 
$$
H_*(B_p,B_{p - 1}) = H_*(B_p, \mathcal{M}^c_p)
$$
Thus by excision we have
$$
H_*(B_p,B_{p - 1}) = H_*(\mathcal{M},\partial  \mathcal{M}_p)
$$
Since any manifold is orientable modulo its boundary with $ \Z_2$ coefficients, we have a nonzero orientation class in $ H_{(n - 1)p + p - 1}(B_p,B_{p - 1})$ which we denote by $ \o_p$.

In contrast with Proposition \ref{p:tt}, we have the following Proposition:

\begin{pro}\label{p:ta}
Under the assumption that (P) has no solution , we have ,
$$
 \mbox{for every } \, p \in \N^* \quad (f_p(\l))_*(\o_p) \not\equiv 0  .
$$
\end{pro}

\begin{pf}

An abstract topological argument displayed in \cite{bc} , pp 260-265 , see also \cite{bb}, which extends virtually to our framework shows that:

$$
\mbox{If } \, (f_1(\l)_* \not\equiv 0 \quad \mbox{then } (f_p(\l))_* \not\equiv 0 \, \mbox{for every } p \geq 2.
$$

Since $ J_{S + \e}$ , for $ \e > 0$ small enough satisfies $ J_{S + \e} \subset V(1,\d)$ , where $ \d \to 0$ if $ \e \to 0$ , one can define using Lemma \ref{l:rep} a continuous map $ s : J_{S + \e} \to \partial M$ which associates to $  u = \ov{\a} \var_{\ov{a},\ov{\l}} + v \in J_{S + \e} \to a \in \partial M $.
Here $ (\ov{\a},\ov{a},\ov{\l})$ are the unique solution of the minimization : $ \min \{ \barre{u - \a \var_{a,\l}} , \, \a \geq 0, \l>0, a \in \partial M \}$

So if $ r: W_1 \to J_{S + \e}$ denotes the retraction by deformation of $ W_1$ onto $ J_{S + \e}$, the existence of a such retraction by deformation follows from the assumption that (P) has no solution from one part and from Proposition \ref{p:ps} from another part.
Let us observe that $ s \circ  r \circ f_1(\l) = id_{\partial M}$
hence $ (f_1(\l))_*(\o_1) \not\equiv 0$ , where $ \o_1$ is the orientation class of $ \partial M$.
Therefore the proof of Proposition \ref{p:ta} is reduced to the abstract topological argument of Bahri-Coron \cite{bc}.
\end{pf}

\vskip .2truecm
\noindent \begin{pfn}{  \sc{of Theorem \ref{t:mt} completed}}

Proposition \ref{p:ta} is in contradiction with Proposition \ref{p:tt}. Therefore (P) has a solution and Theorem \ref{t:mt} is thereby established.
\end{pfn}

\section{Appendix}

\begin{lem}\label{a1}
There holds
$$
\frac{1}{\l^{\frac{n - 2}{2}}} \, \int_{\partial M \cap B_{\rho}(a_1)} \, \d_{a_1,\l}^{\frac{2}{n - 2}}\, \d_{a_2, \l} dv_{g_0} \, = \,o( \int_{\partial M } \, \d_{a_1,\l}^{\frac{n}{n - 2}}\, \d_{a_2, \l} dv_{g_0})  \,   
$$

\end{lem}

\begin{pf}
For $ \e > 0$ a fixed number, let 
$$
 A_{\e} = \{ x \in B_{\rho}(a_1) \cap \partial M ; \d_{a_1,\l} \geq
 \frac{1}{\e \l^{\frac{n - 2}{2}}} \}
$$
Then

\begin{eqnarray*}
\frac{1}{\l^{\frac{n - 2}{2}}} \, \int_{\partial M \cap B_{\rho}(a_1)} \, \d_{a_1,\l}^{\frac{2}{n - 2}}\, \d_{a_2, \l} dv_{g_0} \, 
&\leq & \e \, \int_{A_{\e}} \, \d_{a_1,\l}^{\frac{n}{n - 2}}\, \d_{a_2, \l} dv_{g_0} \,  +\frac{1}{\l^{\frac{n - 2}{2}}} \, \int_{\partial M \cap B_{\rho}(a_1) \setminus A_{\e}} \, \frac{1}{\l \e^{\frac{2}{n - 2}}}\, \d_{a_2, \l} dv_{g_0} \, \\
& \leq & \e  \, \int_{A_{\e}} \, \d_{a_1,\l}^{\frac{n}{n - 2}}\, \d_{a_2, \l} dv_{g_0} +\frac{1}{\l^{\frac{n - 2}{2}}} \, \int_{\partial M \cap B_{\rho}(a_1) \setminus A_{\e}} \, \frac{1}{\l \e^{\frac{2}{n - 2}}}\, \d_{a_2, \l} dv_{g_0} \, \\
& \leq & \, \e \,\int_{\partial M} \, \d_{a_1,\l}^{\frac{n}{n - 2}}\, \d_{a_2, \l} dv_{g_0}  + \, \frac{1}{\e^{\frac{2}{n - 2}}} O(\frac{1}{\l^{\frac{n}{2} + \frac{(n - 2)^2}{2 n}}})
\end{eqnarray*}
Since $\frac{n}{2} + \frac{(n - 2)^2}{2 n} > n - 2 $ and since 
$$
\int_{\partial M} \, \d_{a_1,\l}^{\frac{n}{n - 2}}\, \d_{a_2, \l} dv_{g_0}  \geq \frac{C}{\l^{n - 2}}
$$
then our Lemma follows.
\end{pf}

\begin{lem}\label{l:ineq}[\cite{ba}]
Let $ q > 2$ be given. There exists $ \gamma > 1$ such that  for any $ (a_1, \cdots,a_p) , a_i > 0$, we have
$$
\left( \sum_{i = 1}^p \a_i \right)^q \, \geq \, \sum_{i = 1}^p \a_i^p \, + \, \frac{\g q}{2} \, \sum_{i \not= j} a_i^{q - 1} a_j
$$

\end{lem}

\begin{lem}\label{l:cp} \cite{E2} [Maximum Principle]

Under the assumption that $R_g \geq 0 , \, h_g \geq 0 $ and $ R_g > 0, \, \mbox{or } h_g > 0$ , let $ u \in C^2(\mathring{M}) \cap C^1(\partial M)$ satisfying
$$
L_g u \geq 0 \quad \mbox{on } \mathring{M} \, \quad \mbox{and } B_g u \leq 0 \quad \mbox{on } \partial M
$$
Then $ u \leq 0 \quad \mbox{on } M$

\end{lem}

\par
\par
\par
\par
\par\par\par
\par
\par
\par
\par
\bigskip\bigskip

\noindent
{\bf{Mohameden Ould Ahmedou}} :\\
 Rheinische Frierich-Wilhelms- Universit\"{a}t Bonn \\
Mathematisches Institut, Beringstrasse 4, D- 53115 Bonn, Germany. \\
E-mail: \texttt{ahmedou@math.uni-bonn.de}

\end{document}